\documentclass{article}

\usepackage{amssymb}
\usepackage{amsfonts,amsthm}
\usepackage[dvips]{graphicx}
\usepackage{mathrsfs}
\usepackage{color}

 \newtheorem{theorem}{Theorem}
 \newtheorem{proposition}{Proposition}
 \newtheorem{lemma}{Lemma}
 
 \newtheorem{definition}{Definition}

\newtheorem{assumption}{Assumption}

\newfont{\Bbe}{msbm9}
\newfont{\Bb}{msbm10 scaled\magstep{1}}

 \newcommand{\Z}{\mbox{\Bb Z}}
 \newcommand{\R}{\mbox{\Bb R}}
 \newcommand{\C}{\mbox{\Bb C}}

 \newcommand{\cO}{{\cal O}}
 \newcommand{\cD}{{\cal D}}
 
 \newcommand{\cE}{{\cal E}}
 
 \newcommand{\cU}{{\cal U}}
  \newcommand{\cG}{{\cal G}}
 \newcommand{\cH}{{\cal H}}

\def\phib{{\mbox{\boldmath$\phi$}}}
 \newcommand{\bF}{{\bf F}}
 \newcommand{\bB}{{\bf B}}
 \newcommand{\bD}{{\bf D}}
 \newcommand{\bH}{{\bf H}}
\newcommand{\bE}{{\bf E}}
\newcommand{\bG}{{\bf G}}

 \newcommand{\bU}{{\bf U}}
 
 \newcommand{\bz}{{\bf 0}}
\newcommand{\bn}{{\bf n}}
\newcommand{\bbw}{{\bf w}}

\newcommand{\darr}[4]{{\left\{\begin{array}{ll}
   {#1}&{#2}\\
   {#3}&{#4}
 \end{array}\right.}}

 \newcommand{\ep}{{\epsilon}}

 \newcommand{\diver}{{\rm div}}
 \newcommand{\curl}{{\bf curl\,}}
 \newcommand{\grad}{{\bf grad\,}}
 
 \newcommand{\inprod}[2]{{\langle{#1},{#2}\rangle}}
 \newcommand{\parder}[2]{{\frac{\partial{#1}}{\partial{#2}}}}

 \newcommand{\half}{\frac{1}{2}}
 
 \newcommand{\ia}{({\rm i})}
 \newcommand{\ib}{({\rm ii})}

\newcommand{\hm}{{(H^1_0(D))^m}}

\newcommand{\als}{{ {\rm \; a.s. }}}

\newcommand{\bu}{{\bf u}}
\newcommand{\bbf}{{\bf f}}
\newcommand{\bw}{{\bf w}}
\newcommand{\bv}{{\bf v}}
\newcommand{\by}{{\bf y}}
\newcommand{\bp}{{\bf p}}
\newcommand{\bq}{{\bf q}}
\newcommand{\br}{{\bf r}}
\newcommand{\be}{{\bf e}}

\newcommand\ua{u_{\alpha}}

\newcommand\Zd{{\mathbb Z}^{d}}
\newcommand\E[1]{{\mathbb E}\left[ #1 \right]}
\newcommand{\wsc}{\stackrel{*}{\rightharpoonup}}

\begin{document}
\title{Homogenization of random elliptic systems with an
application to Maxwell's equations\footnote{Research supported by the programme THALES, Code
number 3570 (``Analysis, Modeling and Simulations of Complex and
Stochastic Systems'', University of Crete - Scientific coordinator: M.
Katsoulakis) funded by the European Social Fund and National Hellenic
Sources.} }
\author{ G. Barbatis$^{(1)}$, I. G. Stratis$^{(1)}$ and A. N.
Yannacopoulos$^{(2)}$\\ \\ {\normalsize (1)  Department of Mathematics, National and
Kapodistrian University of Athens, Greece} \\ {\normalsize (2) Department of Statistics,
Athens University of Economics and Business, Greece} }
\date{}
\maketitle

\

\begin{center}

{\em Dedicated to Professor Christodoulos Athanasiadis on the occasion of his retirement}

\end{center}

\

{\bf Abstract:} {\em We study the homogenization of elliptic systems of equations in divergence form where the coefficients are compositions of periodic functions with a random diffeomorphism with stationary gradient.
This is done in the spirit of scalar stochastic homogenization by Blanc, Le Bris and P.-L. Lions. An application of the abstract result is given for Maxwell's equations in random dissipative bianisotropic media.   } 
{\bf MSC 2010 (AMS):  35B27, 35Jxx, 35R60, 60H25, 60H30, 78A48, 78M40}

\section{Introduction}

In a variety of applications, ranging from mechanics to financial engineering,
the mathematical models which arise are in the form of partial differential equations
with variable coefficients which have either a fast periodic variation, or
quasiperiodic variation or even random variation. Relevant examples may be
heat flow is periodic media with fast varying microstructure, pricing of
contingent claims whose underlyings exhibit multiscale volatility or  the modeling of
electromagnetic fields in random complex media. Since partial differential equations with
fast or random variable coefficients present considerable difficulties in their analytic and even
numerical treatment, an approximate model exhibiting constant coefficients would be highly desirable.
This is the main focus of homogenization theory, which
can mathematically be described as a theory for averaging
partial differential equations, whose physical principles have a long history going back to Poisson, Mossotti, Maxwell, Clausius and Rayleigh.

It
is interesting to note that (mathematical) homogenization theory started in the
late 1960s and
in at least three directions, as pointed out by Allaire (see \cite{Allaire} where a detailed reference list can be found):
the first (and oldest) direction dealing with a general theory for the
convergence of operators (namely the H- or G-convergence), the second  treating the asymptotic study of
perforated domains containing many small holes and the third
 refers to a systematic study of periodic structures using
asymptotic analysis, while since the mid 1970's
there is also a variational theory of homogenization, known as
$\Gamma$-convergence.  Mathematical homogenization theory is still a very active field, with a huge number of publications by many researchers  in various directions.

The initial developments in the theory were concerned with homogenization of periodic
structures using either asymptotic analysis or variational tools, but soon the theory was extended to the study of
quasiperiodic structures, adding more realistic features to the models. An important branch of
homogenization theory was concerned with the modeling of random media, in which ergodic theory
plays an important role. Random structures naturally appear in a number of applications (e.g. in the
modeling of porous media or structures with irregular imperfections) and the construction of an
``effective'' homogeneous medium which serves as an approximation of the properties of the full medium
can often be quite useful. What often comes as a surprise is that for a variety of random media, the average medium
is deterministic, a property that arises naturally from the use of ergodic theory.

Historically, stochastic homogenization  for elliptic equations was first studied by Papanicolaou and Varadhan \cite{pava},  and by Kozlov \cite{koz1}.  An important development in random homogenization
theory was the work of Blanc, Le Bris  and Lions \cite{Lions2,Lions} which essentially combined random homogenization and periodic homogenization
for scalar linear elliptic problems by assuming random media whose diffusion coefficient is the composition of
a periodic function with a random diffeomorphism. See, also, \cite{gloria}, \cite{LeTh}. This creates random media which in some sense are ``small'' random
perturbations of periodic structures, thus allowing the extension of the powerful tools of the periodic theory in conjunction
with those of ergodic theory to obtain detailed information concerning the homogenized medium. There is related interesting work on integral  functionals as well as on discrete linear elliptic equations; for both see the recent Habilitation Thesis of Gloria \cite{gloria} and references therein. Quantitative results are also obtained for linear elliptic equations, where connection is made between the statistical properties of the random medium (such as correlation length)  with the properties of the homogenized medium. The passage from discrete to continuous relies basically on the De Giorgi-Nash-Moser theory. Questions related to the convergence rate have been studied both in the linear case, again see the references in \cite{gloria}, and in a very general nonlinear setting by Caffarelli and Souganidis, \cite{casou}.
A different very interesting view point regarding diffusion in random media can be found in \cite{P}.

A number of interesting models require the study of random elliptic systems, rather than scalar equations. As examples one
may consider applications in mechanics (e.g., elasticity, elastoplasticity, thermoviscoplasticity) or in electromagnetic theory\footnote{Homogenization
in electromagnetics has been around for 130 years, which is a testament in
itself to the success and usefulness of this research area.  One of the
first applications is due to James Clerk Maxwell himself (\cite{Max}, p.
401).}  and in particular in the modeling of
complex electromagnetic media, an area with important applications in modern technology (e.g., metamaterials, split-ring arrays, current driven homogenization).  The theory of elliptic systems is rather more involved
that that of scalar equations and it is the aim of the present paper to extend the theory of homogenization of random elliptic equations to that
of elliptic systems. We take as our starting point the modeling assumptions employed by Blanc,  Le Bris and Lions \cite{Lions2,Lions}  i.e. that the diffusivity tensor
is the composition of a periodic tensor with a random diffeomorphism with stationary gradient, an assumption that allows the use of ergodic theory
as a means of obtaining asymptotic results. We show that an extension of homogenization theory for elliptic systems is possible, providing explicit
expressions for the homogenized medium in terms of solutions of properly selected ``cell'' systems, not defined as in the case of classical deterministic homogenization on a periodic cell but rather on the whole space. The results are valid for elliptic systems in general form
and therefore may be used in a wide range of applied models. As an illustration we present a concrete application of the general theory in the
study of homogenization of random bianisotropic media, modeled by the Maxwell equations complemented with constitutive laws in the form of
temporal convolutions involving kernels with random spatial dependence. Such models are often used to study dispersive complex non-homogeneous media, exhibiting random  imperfections.
It is shown that the original homogenization problem can be solved in terms of an auxiliary homogenization problem involving an elliptic system, and using the
general theoretical framework developed in the first part of the present work, we provide expressions for the homogenized medium.

\section{Homogenization of random elliptic systems: Abstract results}
\label{sec:ell}

\subsection{A model for random media}
In this section we present a general class of models for a random medium. Let us consider a general $m\times m$ elliptic system on a bounded  domain $\cO\subset\R^d$, of the form
\begin{equation}
\label{system-0}
\darr{ -\sum_{\alpha=1}^m\sum_{i,j=1}^d \parder{}{x_j}\Big\{ a_{ij\alpha\beta}(x)\parder{u_{\alpha}}{x_i}\Big\}=
f_{\beta},}{\mbox{ in }\cO,}
{u_{\beta}=0,}{\mbox{ on }\partial \cO,}\quad\beta=1,\ldots,m.
\end{equation}
where $u : \cO \to \R^{m}$,  $u=(\ua)$, $\alpha=1,\cdots, m$, is a vector field describing the state of the system, and $A:=(a_{ij\alpha\beta}(x))$, $i,j=1,\cdots, d$, $\alpha,\beta=1,\cdots,m$, is the diffusivity tensor.
Systems of equations of this type appear in a great variety of applications, ranging from elasticity,  fluid mechanics, plasticity, to complex media electromagnetics. The medium (which is considered here in a very general manner) is fully characterized by the diffusivity tensor $A$. In some cases the tensor $A$ cannot be fully described, either on account of unknown processes taking place within the medium or on account of incomplete information. Such cases are considered under the general terminology of random media, and are modelled by assuming that the diffusivity tensor $A$ is a family of random variables; a random field $A(\cdot , \omega):=( a_{ij\alpha\beta}(\cdot , \omega))$, $i,j=1,\cdots, d$, $\alpha,\beta=1,\cdots,m$ where $(\Omega, {\cal F},P)$ is a probability space. The probability space $(\Omega, {\cal F},P)$ is a model for the spatial randomness of the medium; we will return to specific examples later on. The random nature of the medium is used to model imperfections of the medium due to its construction, experimental mis-specifications of the medium properties etc. Clearly, if $A$ is a random field, then the solution of system (\ref{system-0}) is a random field as well.  The essential meaning of this formulation is that for a random medium, if an experiment is repeated, then, we expect different results as outcome of the experiment. Then, in $\Omega$ we collect all the possible outcomes that an experiment concerning this medium may provide, and by $\omega \in\Omega$ we denote a particular outcome of an experiment.

We now restrict the above very general setup to one which is more useful for the study of homogenization in random media. We wish to restrict our study to random media that present some sort of self-repetitive structure, i.e., of some structure that allows us to reconstruct (in a statistical manner of speaking) the whole medium from knowledge of a part of it only. This certainly is true for media presenting some periodic structure, in which case the whole medium can be reconstructed (exactly) by knowledge of a ``fundamental cell'' and by translations of it by an appropriate vector. Clearly, this cannot be true for a random medium, however, we may provide a convenient framework that allows us to do that in an approximate (statistical) fashion. This can be done within the framework of ergodic or stationary media.

Let $(\Omega, {\cal F},P)$ be a probability space, and $\tau_{k} : \Omega \to \Omega$ a group of transformations parameterized by a parameter $k \in \Z^d$. We assume that the group of transformations $\{\tau_{k}\}$   preserves the measure $P$, i.e.,
\begin{eqnarray}
P(\tau_{k} A) = P(A), \,\,\,\,\ \forall A \in {\cal F},  \,\,\, \forall k \in \Z^d.
\nonumber
\end{eqnarray}
The probability space $(\Omega, {\cal F},{\mathbb P})$ is to be interpreted as follows: Each realization $\omega$ is to be interpreted as a particular configuration of the medium. In other words, each experiment we perform on a particular medium corresponds to a particular choice of $\omega \in \Omega$. However, it is not known beforehand and with certainty which medium is to be realized, when the experiment is performed. The probability that a particular medium is realized is given by the probability measure ${\mathbb P}$. The above description is rather abstract, and accomodates a number of interesting  cases arising in applications. A more concrete description is to assume that
$\Omega=\R^{d}$ i.e. each $\omega$ is identified with a point $x \in {\mathbb R}^{d}$, and assume that $\tau$ is  $({\mathbb R}^{3},+)$ the usual translation group, $\tau_{k}\omega=\tau_{k}x=x+k$ the action of a group acting on $\Omega$. Invariance of measure means some kind of periodicity with respect to a fundamental lattice, i.e., self repetitive structure, obtained by translation by $k \in \Zd$ of a fundamental structure.

In this paper we will assume certain conditions on the random coefficients. These are  the conditions of ergodicity and stationarity.

\begin{definition}[Stationarity and ergodicity] \hfill
\begin{enumerate}

\item The group action $(\tau_k)$ is ergodic if it is measure preserving and any invariant event $A$ has probability $0$ or $1$, i.e.
\begin{eqnarray}
(  \tau_k \,A=A  \; , \;\;  \forall k \in \Z^d) \,\, \Longrightarrow \,\, (P(A)=0  \mbox{ or } P(A)=1)
\nonumber
\end{eqnarray}

\item A random field $F \in L^{1}_{loc}(\R^{d},L^{1}(\Omega))$ is called stationary with respect to the group action $(\tau_k)$ if
\begin{eqnarray}
 F(x+k,\omega)=F(x,\tau_{k}\omega), \,\,\;\; \forall k \in \Z^d \mbox{ a.e. in x, a.s.}
\nonumber
\end{eqnarray}

\end{enumerate}
\end{definition}
Stationary processes need not be ergodic (consider for example $F(x,\omega)=Y(\omega)$ where $Y$ is a given random variable).
Stationarity  guarantees that in a statistical sense, parts of the material located at different positions will present the same properties, i.e. that the statistical properties of the medium are invariant under translations which are to be understood as the transformation $\tau_{x}$. In other words, the function $F$ at $x \in \R^{d}$ and the function $F$ at $x+y$, $y \in \R^d$,  will look as if generated by the same probabilistic law.  Alternatively, if $F$ is considered as $F: \R^{d} \times \Omega \to \R^{N}$ then $F$ is stationary if and only if $P(\{\omega \,: \, F(x, \omega) \in B\})$ is independent of $x$ for $B \in {\cal B}(\R^{N})$. If $f$ is an $L^{1}(\Omega)$ random variable, and we define the random field $F(x,\omega):=f(\tau_{x}\omega)$, then, $F(x,\omega)$ is a stationary random field. Ergodicity implies that all $\tau$ invariant quantities are non-random. For examples of media that fall within this description we refer to \cite{Telega}.

For the purpose of this paper, inspired by recent very interesting  work of Lions, Le Bris and Le Blanc \cite{Lions} on stochastic elliptic homogenization, we will concentrate on random coefficients of a special form.
\begin{assumption}
The coefficients of the problem depend additionally on a small parameter $\epsilon>0$ and are random fields of the form
\begin{eqnarray}
F\left( \frac{x}{\epsilon}, \omega \right)=F_{per}\left( \Phi^{-1}\left(\frac{x}{\epsilon},\omega \right)\right)
\label{eqn:ASSUMPCOEF}
\end{eqnarray}
where $F_{per}$ is a periodic function in ${\mathbb R}^{d}$ and $\Phi$ is a random mapping which is a diffeomorphism almost surely, with stationary gradient.
\end{assumption}
This type of coefficients models some sort of statistical periodicity of the medium. Problems of this type have been studied in \cite{Lions2,Lions}.

Stationarity and ergodicity allow us to look at average properties of the material at long scales and obtain nice expressions for these quantities. In fact the Ergodic Theorem (see e.g. \cite{Lions}) states that
\begin{eqnarray*}
 \lim_{N \to \infty} \frac{1}{(2N+1)^{d}} \sum_{|k|_{\infty} \le N} F(x,\tau_{k}\omega) = \E{F(x,\cdot)}, \,\,\, \mbox{in}\,\, L^{\infty}(\R^{d}), \,\,\, \mbox{a.s.},
 \end{eqnarray*}
 where for $k=(k_{i}) \in \Zd$, $i=1,\cdots, d$, $|k|_{\infty}=\max_{i} |k_{i}|$.  This implies that
\begin{eqnarray}
F\left( \frac{x}{\epsilon},\omega \right) \wsc
\E{\int_{Q}F(x, \cdot) dx}, \,\,\, \mbox{as} \,\, \epsilon \rightarrow 0\,\, \mbox{a.s. \, in} \,\, L^{\infty}({\mathbb R}^{d})
\nonumber
\end{eqnarray}
where $\E{\cdot}$ is the expectation over the measure $P$, $\wsc$ denotes the weak star convergence and $Q$ is the unit cube in $\R^{d}$.
The ergodicity hypothesis implies that instead of looking at an ensemble average of media, and averaging the properties of the medium on the ensemble average, we may consider a single realization of the medium whose spatial dimensions are large  and sample its properties   by traversing this single realization for large enough distances.

\subsection{Homogenization of the random elliptic system}\label{sec:ELLHOM}

We consider the following $m\times m$ system in a bounded Lipschitz domain $\cO\subset\R^d$:
\begin{equation}
\label{system}
\darr{ -\sum_{\alpha=1}^m\sum_{i,j=1}^d \parder{}{x_j}\Big\{ a_{ij\alpha\beta}(x)\parder{u_{\alpha}}{x_i}\Big\}=
f_{\beta},}{\mbox{ in }\cO,}
{u_{\beta}=0,}{\mbox{ on }\partial \cO,}\quad\beta=1,\ldots,m.
\end{equation}
The functions $a_{ij\alpha\beta}$ are deterministic functions in $L^{\infty}(\cO)$.
The system is understood in the weak sense: let $Q(\cdot,\cdot)$ be the bilinear form
on $(H^1_0(\cO))^m$ given by
\[
Q(\bu,\bv)=\int_{\cO}\sum_{\alpha,\beta=1}^m\sum_{i,j=1}^d a_{ij\alpha\beta}(x)\parder{u_{\alpha}}{x_i}
\parder{v_{\beta}}{x_j}dx \, , \;\;\; \bu=(u_{\alpha}),\bv=(v_{\beta})\in\hm\, ;
\]
to simplify our notation we shall also write this as
\begin{equation}
Q(\bu,\bv)=\int_{\cO} A(x)\nabla\bu\cdot \nabla\bv \, dx\, .
\label{a}
\end{equation}
Thus, given $\bbf=(f_{\beta})_{\beta=1}^m \in (H^{-1}(\cO))^m$ we say that $\bu=(u_{\alpha})\in (H^1_0(\cO))^m$ is a (weak) solution of the system
(\ref{system}) if
\[
Q(\bu,\bv)=\inprod{\bbf}{\bv}:=\sum_{\beta=1}^m \inprod{f_{\beta}}{v_{\beta}}
\]
for all $\bv=(v_{\beta})\in (H^1_0(\cO))^m$. We make the ellipticity assumption that there exists $c>0$ such that
\[
Q(\bu,\bu)\geq c\| \bu\|_{\hm}^2 , \;\; \;\; \bu\in\hm\, ;
\]
hence the Lax-Milgram theorem yields the existence of a unique solution  $\bu$ of (\ref{system}), which satisfies
$\|\bu\|_{\hm}\leq c^{-1}\|\bbf\|_{(H^{-1}(\cO))^m}$.

{\em Notational remark.} We think of $A(x)=(a_{ij\alpha\beta}(x))$ as a matrix acting on $\R^{md}$: if $\bp=(p_{i\alpha})\in\R^{md}$, then
$\bq=A\bp$ is the vector $q_{j\beta}=a_{ij\alpha\beta}p_{i\alpha}$. Here and below, we use the summation convention over repeated indices.
Moreover, Latin letters $i,j,k\ldots$ will take values in $\{1,\ldots,d\}$, while Greek letters $\alpha,\beta,\gamma,\ldots$ will take values in
$\{1,\ldots,m\}$.

Our aim in this section is to study the homogenization problem for the system (\ref{system}) when the coefficients
are random of the form (\ref{eqn:ASSUMPCOEF}). More precisely, for $\epsilon>0$ we consider the stochastic elliptic system
\begin{equation}
\label{syst_eps}
\darr{ -\sum_{\alpha=1}^m\sum_{i,j=1}^d \parder{}{x_j}\Big\{ a_{ij\alpha\beta}(\Phi^{-1}(\frac{x}{\epsilon},\omega))\parder{u_{\alpha}}{x_i}\Big\}=
f_{\beta},}{\mbox{ in }\cO,}
{u_{\beta}=0,}{\mbox{ on }\partial \cO,}\quad\beta=1,\ldots,m ,
\end{equation}
where the functions $a_{ij\alpha\beta}(y)$ are periodic of period $Q$. We intend
to study the behaviour of the solution $\bu^{\epsilon}(x,\omega)$ of (\ref{syst_eps}) as $\epsilon\to 0$.

In the rest of the paper (weak) convergence in $L^p_{loc}(\R^d)$ means (weak) convergence in $L^p(U)$ for any $U\subset\subset\R^d$.

\begin{lemma}
\label{lem:1}
We have as $\epsilon\to 0$,
\begin{eqnarray*}
\ia && (\nabla\Phi)(\frac{x}{\epsilon},\omega)\longrightarrow \E{\int_Q(\nabla\Phi)(y,\cdot)dy}\quad
\mbox{\rm $*$-weakly in }(L^{\infty}(\R^d))^{d\times d}\, , \; \als \, ; \\
\ib && \epsilon \Phi(\frac{x}{\epsilon},\omega)\longrightarrow \E{\int_Q(\nabla\Phi)(y,\cdot)dy}x\quad
\mbox{\rm in }(L^{\infty}_{loc}(\R^d))^d\, , \; \als \, .
\end{eqnarray*}
\end{lemma}
{\em Proof.} Part $\ia$ is an immediate consequence of the Ergodic Theorem. To prove $\ib$, let us define
$\Psi_{\epsilon}(x,\omega)=\epsilon \Phi(\frac{x}{\epsilon},\omega)- \E{\int_Q(\nabla\Phi)(y,\cdot)dy}x$.
Let $U\subset\subset\R^d$ be fixed. Then $\|\Psi_{\epsilon}\|_{L^{\infty}(U)}\leq c$ and, by $\ia$,
$\nabla \Psi_{\epsilon}\rightharpoonup 0$ $*$-weakly in $(L^{\infty}(\R^d))^{d\times d}$. Hence
\[
\|\Psi_{\epsilon}\|_{W^{1,\infty}(U)}\leq c \, , \;\; \epsilon>0\, , \; \als\, .
\]
Therefore there exists $\Psi\in W^{1,\infty}(U)$ such that, up to a subsequence, $\Psi_{\epsilon} \rightharpoonup \Psi$
weakly in $W^{1,\infty}(U)$ and hence $\Psi_{\epsilon}\to\Psi $ (strongly) in $L^{\infty}(U)$. It then follows immediately that $\Psi=0$ and
that the full sequence $\Psi_{\epsilon}\to 0$. $\hfill \Box$

\begin{lemma}
\label{lem:n1}
Let $\psi$, $\psi_{\epsilon}$, $\epsilon>0$, be diffeomorphisms of $\R^d$ onto itself such
that
\[
\max\{\|\nabla\psi_{\epsilon}\|_{\infty},\|\nabla\psi_{\epsilon}^{-1}\|_{\infty}\}
\leq c_1
\]
for all $\epsilon>0$ and $\psi_{\epsilon}\to \psi$ in $L^{\infty}_{loc}(\R^d)$. Then
\[
\chi_{\psi_{\epsilon}(A)} \longrightarrow \chi_{\psi(A)}\, , \;\; \mbox{ in }L^1(\R^d),
\]
for any set $A\subset\subset\R^d$.
\end{lemma}
{\em Proof.} Let us denote by $B^{\delta}$ the $\delta$-neighbourhood of a set $B\subset\R^d$. Let $A\subset\subset\R^d$ be given and let
$\delta(\epsilon)=\|\psi_{\epsilon}-\psi\|_{L^{\infty}(A)}$. Then
\[
{\rm Hausdorff \, dist}(\psi_{\epsilon}(A),\psi(A)) \leq\delta(\epsilon),
\]
and therefore
\[
\|\psi_{\epsilon}(A) -\psi(A)\|_{L^1(\R^d)}=
| \psi_{\epsilon}(A)\bigtriangleup \psi(A)| \leq | (\psi_{\epsilon}(A))^{\delta(\epsilon)}\setminus
\psi_{\epsilon}(A)| + | \psi(A)^{\delta(\epsilon)}\setminus \psi(A)| .
\]
The second of the last two terms clearly tends to zero as $\epsilon\to 0$. Moreover, it is easily seen that
$\psi_{\epsilon}(A)^{\delta}\subset \psi_{\epsilon}(A^{c_1\delta})$ for any $\delta>0$. Hence
\[
| \psi_{\epsilon}(A)^{\delta(\epsilon)}\setminus \psi_{\epsilon}(A)|
\leq | \psi_{\epsilon}(A^{c_1\delta(\epsilon)})\setminus \psi_{\epsilon}(A)|\leq
| (\psi_{\epsilon}(A^{c_1\delta(\epsilon)}\setminus A)|\leq c_1 |A^{c_1\delta(\epsilon)}\setminus A|\to 0,
\]
as $\epsilon\to 0$. This concludes the proof. $\hfill\Box$

We set
\[
c_{\Phi}=\det\bigg( \E{\int_Q \nabla\Phi(x,\cdot)dx}\bigg)^{-1}.
\]
\begin{lemma}
Let $g\in L^{\infty}(\R^d, L^1(\Omega))$ be stationary. Then
\[
g(\Phi^{-1}(\frac{x}{\epsilon},\omega),\omega) \rightharpoonup c_{\Phi}\E{\int_{\Phi(Q,\cdot)}g(\Phi^{-1}(x,\cdot),\cdot)dx} ,
\]
as $\epsilon\to 0$, $*$-weakly in $L^{\infty}(\R^d)$, almost surely.
\label{lem:2}
\end{lemma}
{\em Proof.} A simple approximation argument shows that it suffices to check the convergence against functions of
the form $\chi_A$, where $A\subset\R^d$ is open, bounded and connected. Now, Lemma \ref{lem:1} (ii) and Lemma \ref{lem:n1}
imply that for such a set $A$ we have
\begin{equation}
\chi_{\epsilon\Phi^{-1}(\frac{A}{\epsilon},\omega)} \to \chi_B  \, , \;\; \mbox{ in }L^1(\R^d),
\label{3}
\end{equation}
as $\epsilon\to 0$, almost surely, where
\[
B=\Big( \E{ \int_Q\nabla\Phi(y,\cdot)dy}\Big)^{-1}A\, .
\]
Moreover, the Ergodic Theorem applied to the stationary random function $F(x,\omega)=g(x,\omega)\det\nabla\Phi(x,\omega)$ gives
\begin{eqnarray}
g(\frac{y}{\epsilon},\omega)\det\nabla\Phi(\frac{y}{\epsilon},\omega) &\wsc& \E{\int_Q g(y,\cdot)\det\nabla\Phi(y,\cdot)dy} \nonumber\\
&=& \E{\int_{\Phi(Q,\cdot)} g(\Phi^{-1}(x,\cdot),\cdot)dx},
\label{4}
\end{eqnarray}
as $\epsilon\to 0$, $*$-weakly in $L^{\infty}(\R^d)$, almost surely. Combining (\ref{3}) and (\ref{4}) we obtain
\begin{eqnarray*}
&& \hspace{-2cm}\int_A g(\Phi^{-1}(\frac{x}{\epsilon},\omega),\omega)dx =\\
&=&\int_{\R^d}g(\frac{y}{\epsilon},\omega)\det\nabla\Phi(\frac{y}{\epsilon},\omega)\chi_{\epsilon\Phi^{-1}(\frac{A}{\epsilon},\omega)}dy\\
&\to& \int_{\R^d}  \E{\int_{\Phi(Q,\cdot)} g(\Phi^{-1}(x,\cdot),\cdot)dx}\chi_B(y)dy\\
&=&  \E{\int_{\Phi(Q,\cdot)} g(\Phi^{-1}(x,\cdot),\cdot)dx}|B|\\
&=&  \E{\int_{\Phi(Q,\cdot)} g(\Phi^{-1}(x,\cdot),\cdot)dx}
\det\bigg( \E{ \int_Q \nabla\Phi(x,\cdot)dx }\bigg)^{-1}|A|,
\end{eqnarray*}
as required. $\hfill\Box$

We now need to consider a corrector problem which shall play a crucial role for the main theorem of this section. We note that in the periodic case the corrector problem (\ref{aux}) below is posed on a single cell, but in the general stationary case it is posed globally on $\R^d$.

For $0<\delta \leq 1$ we denote by $C^{\delta}(\R^d)$ the set of functions $v$ for which the H\"{o}lder seminorm
\[
[v]_{\delta} =\sup_{x,y\in\R^d , x\neq y} \frac{ |u(y)-u(x)|}{|y-x|^{\delta}}
\]
is finite.
\begin{proposition}
Assume that the functions $a_{ij\alpha\beta}$ belong in $C^{\delta}(\R^d)$ for some $0<\delta\leq 1$.  Then for any $\bp=(p_{i\alpha})\in\R^{md}$, the system
\begin{equation}
\label{aux}
\left\{
\begin{array}{l}
{ -\sum_{\alpha=1}^m\sum_{i,j=1}^d \parder{}{y_j}\Big\{ a_{ij\alpha\beta}(\Phi^{-1}(y,\omega))(\parder{w_{\alpha}}{y_i}+p_{i\alpha})
\Big\}=0,\;\mbox{ in }\R^d,}{\;\; (\beta=1,\ldots,m) }\\[0.2cm]
{\bw(y,\omega)=\tilde{\bw}(\Phi^{-1}(y,\omega),\omega),\;\;\; \nabla \tilde{\bw}
\mbox{ stationary, }}{}\\[0.2cm]
{\E{\int_{\Phi(Q,\cdot)}\nabla \bw(y,\cdot)dy }=0,}{}
\end{array}
\right.
\end{equation}
has a unique (up to an additive constant which may depend on $\omega$) solution $\bw^{(\bp)}$ in $(H^1_{loc}(\R^d, L^2(\Omega)))^m$.
\label{prop:1}
\end{proposition}
{\em Proof.} {\em Existence.} We regularize the problem (\ref{aux}) by fixing $\theta >0$ (which will eventually tend to zero) and considering the system
\begin{equation}
\left\{
\begin{array}{l}
{-\sum_{\alpha=1}^m\sum_{i,j=1}^d \parder{}{y_j}\Big\{ a_{ij\alpha\beta}(\Phi^{-1}(y,\omega))(\parder{w_{\alpha}}{y_i}+p_{i\alpha})\Big\}+\theta
 w_{\beta}=0 ,}\\[0.2cm]
{\bw(y,\omega)=\tilde \bw(\Phi^{-1}(y,\omega),\omega),}{\quad \tilde \bw \mbox{ stationary.} }
\end{array}\right.
\label{theta}
\end{equation}
The problem (\ref{theta}) is understood as follows: we define the Hilbert space
\[
\cH=\{ \bw=\tilde \bw\circ\Phi^{-1}\; : \; \tilde \bw\in H^1_{loc}(\R^d, L^2(\Omega)) \, , \; \mbox{ $\tilde\bw$ stationary}\},
\]
equipped with the (real) inner product
\[
\inprod{\bw}{\bw'}_{\cH}=\E{ \int_Q  \sum_{\alpha=1}^m\Big(  \sum_{i=1}^d \parder{\tilde w_{\alpha}}{y_i}
\parder{\tilde w_{\alpha}'}{y_i} + \tilde w_{\alpha}\tilde w'_{\alpha} \Big)  dy }.
\]
We define the bilinear form
\[
a(\bw,\bw' ;\theta)=\E{\int_{\Phi(Q,\cdot)} \Big(
\sum_{\alpha,\beta=1}^m\Big [\sum_{i,j=1}^d a_{ij\alpha\beta}(\Phi^{-1}(y,\cdot))\parder{ w_{\alpha}}{y_i}
\parder{\tilde w_{\beta}'}{y_j} +\theta \delta_{\alpha\beta} w_{\alpha}\tilde w_{\beta}' \Big] \Big)dy },
\]
and the problem (\ref{theta}) is to find $\bw\in\cH$ so that
\[
a(\bw,\bw' ; \theta) = \E{\int_{\Phi(Q,\cdot)} \sum_{\alpha,\beta=1}^m \sum_{i,j=1}^d a_{ij\alpha\beta}(\Phi^{-1}(y,\cdot))p_{i\alpha}
\parder{w_{\beta}'}{y_j}  dy }
\, , \;\; \forall \bw'\in\cH\, .
\]

Since the form $a(\bw,\bw';\theta)$ is coercive on $\cH$, equation (\ref{theta}) has a unique solution $\bw^{(\bp,\theta)}$. Hence, in the notation used in
(\ref{a}), $\bw^{(\bp,\theta)}$ satisfies
\begin{equation}
\E{\int_{\Phi(Q,\cdot)} \Big( A(\Phi^{-1}(y,\cdot))(\nabla\bw^{(\bp,\theta)} +\bp) \cdot \nabla\bw'  +\theta \bw^{(\bp,\theta)}\cdot\bw'
\Big) dy} = 0\, , \;\; \forall  \bw'\in\cH\, .
\label{wtheta}
\end{equation}
The solution $\bw^{(\bp,\theta)}$ satisfies
\begin{equation}
\E{\int_Q |\nabla\tilde\bw^{(\bp,\theta)}|^2dy}\leq c \; , \quad\quad
\E{\int_Q |\tilde\bw^{(\bp,\theta)}|^2dy}\leq \frac{c}{\theta} \, .
\label{eq:lm}
\end{equation}
In particular
\[
\| \nabla\tilde\bw^{(\bp,\theta)}\|_{L^2(Q, L^2(\Omega))}\leq c\; , \; \;\forall\theta >0,
\]
hence, for a subsequence $\theta\to 0$,
\begin{equation}
\nabla\tilde\bw^{(\bp,\theta)} \rightharpoonup T^{(\bp)}, \qquad \mbox{ weakly in $L^2(Q, L^2(\Omega))$,}
\label{wconv}
\end{equation}
for some $T^{(\bp)}\in L^2(Q, L^2(\Omega))$.
We extend $T^{(\bp)}(y,\omega)$ from $Q$ to $\R^d$ by requiring it to be stationary. Then it is easily verified that
\[
\nabla\tilde\bw^{(\bp,\theta)} \rightharpoonup T^{(\bp)}, \qquad \mbox{ weakly in $L^2_{loc}(\R^d, L^2(\Omega))$,
as $\theta\to 0$.}
\]
Defining
\begin{equation}
\label{5}
S^{(\bp)}(y,\omega) =(\nabla\Phi^{-1})(y,\omega) T^{(\bp)}(\Phi^{-1}(y,\omega),\omega) \, , \; y\in\R^d \, , \; \omega\in\Omega ,
\end{equation}
it is easily seen that, for the same subsequence as above,
\[
\nabla\bw^{(\bp,\theta)} \rightharpoonup S^{(\bp)}, \qquad \mbox{ weakly in $L^2_{loc}(\R^d, L^2(\Omega))$, as $\theta\to 0$.}
\]
{\bf Claim.} There exist $\bv^{(\bp)}$, $\tilde \bv^{(\bp)}$ in $H^1_{loc}(\R^d , L^2(\Omega))$ such that
\[
T^{(\bp)}=\nabla\tilde \bv^{(\bp)} \;\; , \qquad S^{(\bp)}=\nabla \bv^{(\bp)}.
\]
We shall prove the claim for $T^{(\bp)}$, the proof for $S^{(\bp)}$ is similar.
Let us fix $\phi\in C^{\infty}_c(\R^d)$. For each $1\leq i,j\leq d$ and $\theta >0$ we then have
\begin{equation}
\int_{\R^d} {\tilde\bbw}^{(\bp,\theta)}_{x_i} \phi_{x_j}dx =
\int_{\R^d} {\tilde\bbw}^{(\bp,\theta)}_{x_j} \phi_{x_i}dx \; , \quad \omega\in\Omega.
\label{evi}
\end{equation}
We next multiply (\ref{evi}) by an arbitrary $\psi\in L^2(\Omega)$, integrate over $\Omega$ and let $\theta\to 0$ along the subsequence above. Using the weak convergence (\ref{wconv}) we obtain
\[
\E{\int_{\R^d} T^{(\bp)}_i \phi_{x_j}(x)\psi(\omega) dx}  =  \E{\int_{\R^d} T^{(\bp)}_j \phi_{x_i}(x)\psi(\omega) dx} .
\]
Since $\psi\in L^2(\Omega)$ is arbitrary this gives
\begin{equation}
\int_{\R^d} T^{(\bp)}_i \phi_{x_j}dx = \int_{\R^d} T^{(\bp)}_j \phi_{x_i}dx \; , \qquad \mbox{ for all $\omega\in\Omega$ and }\phi\in C^{\infty}_c(\R^d),
\label{dr}
\end{equation}
that is $(T^{(\bp)}_i)_{x_j} = (T^{(\bp)}_j)_{x_i}$ in distributional sense. By \cite[Theorem 2.1]{Mar} this implies that there exists a distribution
$\tilde \bv^{(\bp)}$, which may also depend on $\omega$, such that $T^{(\bp)}=\nabla\tilde \bv^{(\bp)}$. The fact that $ \bv^{(\bp)}\in H^1_{loc}(\R^d , L^2(\Omega))$ follows from \cite[Theorem 3.1]{Mar}.

The Claim together with (\ref{5}) imply
\[
\bv^{(\bp)}(y,\omega)= \tilde \bv^{(\bp)}(\Phi^{-1}(y,\omega),\omega) + X(\omega) \, , \;\; y\in\R^d\, , \; \als\, ,
\]
for some random variable $X(\omega)$, $\omega\in\Omega$, independent of $y\in\R^d$.
Let us now define $\bw^{(\bp)}=\bv^{(\bp)}-X$, $\tilde \bw^{(\bp)}=\tilde \bv^{(\bp)}$. Changing variables in (\ref{wtheta}), taking the limit $\theta\to 0$ and using (\ref{eq:lm}) we obtain
\begin{equation}
\E{\int_Q  A(y)(\nabla \bw^{(\bp)}+\bp) \cdot \nabla\bw'  dy} =  0 \, , \;\; \forall  \bw'\in\cH\, ;
\label{lim}
\end{equation}
Moreover, we have by construction $\bw^{(\bp)}(y,\omega)=\tilde\bw^{(\bp)}(\Phi^{-1}(y,\omega),\omega)$ and
$\E{\int_{\Phi(Q)}\nabla \bw^{(\bp)}(y,\cdot)dy}=0$. Hence existence has been proved.

{\em Uniqueness.} Suppose $\bw$ is a solution of the corresponding homogeneous problem:
\begin{equation}
\label{pa}
\left\{
\begin{array}{l}
{ -\sum_{\alpha=1}^m\sum_{i,j=1}^d \parder{}{y_j}\Big\{ a_{ij\alpha\beta}(\Phi^{-1}(y,\omega))\parder{w_{\alpha}}{y_i}
\Big\}=0,\;\mbox{ in }\R^d,}{\;\; (\beta=1,\ldots,m) }\\[0.2cm]
{\bw(y,\omega)=\tilde{\bw}(\Phi^{-1}(y,\omega),\omega),\;\;\; \nabla \tilde{\bw}
\mbox{ stationary, }}{}\\[0.2cm]
{\E{\int_{\Phi(Q,\cdot)}\nabla \bw(y,\cdot)dy }=0,}{}
\end{array}
\right.
\end{equation}
By standard elliptic regularity, $\nabla\bw\in L^{\infty}_{loc}(\R^d)$ and, therefore, also $\nabla\tilde\bw\in L^{\infty}_{loc}(\R^d)$. We use the fact (see \cite{koz2,AS}) that the stationarity of $\nabla\tilde\bw$ implies that $\tilde\bw$ and $\bw$ are sublinear at infinity,
that is $\bw(y)=o(|y|)$.
Now, let $N$ be a large parameter and $\chi_N$ be a cut-off function which equals 1 on $Q_{N}$, vanishes on $Q_{2N}$ and satisfies $|\nabla\chi_N|\leq c/N$. We multiply the equation in (\ref{pa}) by $\chi_N w_{\alpha}$ and integrate by parts. We obtain for $\omega\in\Omega$
\[
\int_{Q_{2N}}  A(\Phi^{-1}(y,\omega))\nabla\bw  \cdot \nabla (\chi_N\bw )  dy =0.
\]
Using the Cauchy-Schwarz inequality this gives
\[
\int_{Q_{N}}  A(\Phi^{-1}(y,\omega))\nabla\bw  \cdot \nabla\bw  dy \leq\frac{c}{N}\|\nabla \bw\|_{L^2(Q_{2N})} \| \bw\|_{L^2(Q_{2N})}.
\]
Now, by the local regularity estimates of Hong \cite[Theorem 2]{Dong} we have that
$\|\nabla \bw\|_{L^2(Q_{2N})} \leq c N^{-1} \|\bw\|_{L^2(Q_{4N})}$. Combining the above we obtain
\begin{equation}
\label{eid}
\frac{1}{N^d}\int_{Q_{N}}  A(\Phi^{-1}(y,\omega))\nabla\bw  \cdot \nabla\bw  dy \leq cN^{d-2}\| \bw\|_{L^2(Q_{4N})}^2 = o(1) \;  ,
\end{equation}
as $N\to\infty$. We note that the integral in the LHS of (\ref{eid}) equals $\int_{Q_N}g(\Phi^{-1}(y,\omega),\omega)dy$ where
\[
g(z,\omega)=[A(z)(\nabla\Phi(z,\omega))^{-1}\nabla\tilde\bw (z,\omega)](\nabla\Phi(z,\omega))^{-1}\nabla\tilde\bw (z,\omega)
\]
is stationary. Because of stationarity, Lemma \ref{lem:2} applies here. Taking $\epsilon=1/N$ and the test-function $\chi_{Q}\in L^1(\R^d)$ we obtain that
\begin{equation}
\frac{1}{N^d}\int_{Q_N}g(\Phi^{-1}(y,\omega),\omega)dy =\int_{Q}g(\Phi^{-1}(Nx,\omega),\omega)dx \rightarrow
c_{\Phi}\E{\int_{\Phi(Q,\cdot)}g(\Phi^{-1}(x,\cdot),\cdot)dx},
\label{eid1}
\end{equation}
as $N\to \infty$. From (\ref{eid}) and (\ref{eid1}) we conclude that $g=0$ and hence $\bw=0$. This concludes the proof.
$\hfill\Box$

Let $\{\be_{i\alpha}\}$, $\alpha=1,\ldots,m$, $i=1,\ldots,d$, be the canonical basis of $\R^{md}$. When $\bp=\be_{i\alpha}$ we shall write $\bw^{(i\alpha)}$ instead of $\bw^{(\be_{i\alpha})}$.
So each $\bw^{(i\alpha)}$ is a vector field, $\bw^{(i\alpha)}=(w^{(i\alpha)}_{\beta})_{\beta=1}^m$.

We define the homogenized coefficient matrix $A^*=\{a_{ij\alpha\beta}^*\}$ by
\begin{equation}
\label{eq:def_hom_m}
a_{ij\alpha\beta}^* = \det\bigg( \E{ \int_Q \nabla\Phi(x,\cdot)dx }\bigg)^{-1} \E{\int_{\Phi(Q,\cdot)} a_{kj\gamma\beta}\Big(\delta_{ik}
\delta_{\alpha\gamma} + \parder{w^{(i\alpha)}_{\gamma}}{y_{k}} \Big)dy}.
\end{equation}

\begin{theorem}
Let $\cO\subset\R^d$ be bounded and ${\bf f}=(f_{\alpha})_{\alpha=1}^m \in (H^{-1}(\cO))^m$.
Assume that the functions $a_{ij\alpha\beta}$ belong in $C^{\delta}(\R^d)$ for some $0<\delta\leq 1$.
Then the solution $\bu^{(\epsilon)}(x,\omega)$ of (\ref{syst_eps}) converges as $\epsilon\to 0$ weakly in $(H^1_0(\cO))^m$ almost surely to the function $\bu^*(x)$ which is the solution of the problem
\begin{equation}
\darr{ -\sum_{\alpha=1}^m\sum_{i,j=1}^d \parder{}{x_j}\Big\{ a_{ij\alpha\beta}^*\parder{u_{\alpha}}{x_i}\Big\}=
f_{\beta},}{\mbox{ in }\cO,}
{u_{\beta}=0,}{\mbox{ on }\partial \cO,}\quad\beta=1,\ldots,m.
\label{hom_syst}
\end{equation}
\label{thm:elliptic_homog}
\end{theorem}
{\em Proof.}  The sequence $(\bu^{\epsilon})$ is bounded in $L^{\infty}(\Omega,(H^1_0(\cO))^m)$. Hence the sequence
$\br^{\epsilon}(x,\omega):=A(x/\epsilon)\nabla\bu^{\epsilon}$ is bounded in $L^{\infty}(\Omega,(L^2(\cO))^{md})$. It follows that, by extracting a subsequence $\epsilon\to 0$, there exist $\bu^*$ and $\br^*$ such that
\begin{eqnarray}
&& \bu^{\epsilon} \rightharpoonup \bu^* \quad\mbox{ weakly in } (H^1_0(\cO))^m \als \label{p1}\\[0.2cm]
&& \br^{\epsilon} \rightharpoonup \br^* \quad\mbox{ weakly in } (L^2(\cO))^{md} \als . \label{p2}
\end{eqnarray}
The proof of the theorem will be complete once we prove that
\begin{equation}
\br^*= A^* \nabla\bu^* \, , \;\; x\in \cO\, , \; \als \, .
\label{wtp}
\end{equation}
Let $\bp\in \R^{md}$ be fixed and let $\bv^{(\bp)}$ be the solution of the problem
\begin{equation}
\label{auxv}
\left\{
\begin{array}{l}
{ -\sum_{\alpha=1}^m\sum_{i,j=1}^d \parder{}{y_j}\Big\{ a_{ji\beta\alpha}(\Phi^{-1}(y,\omega))(\parder{v_{\alpha}}{y_i}+p_{i\alpha})
\Big\}=0,\;\mbox{ in }\R^d,}{ \;\; (\beta=1,\ldots,m)}\\[0.2cm]
{\bv(y,\omega)=\tilde{\bv}(\Phi^{-1}(y,\omega),\omega),\;\;\; \nabla\tilde{\bv} \mbox{ stationary, }}{}\\[0.2cm]
{\E{\int_{\Phi(Q,\cdot)}\nabla \bv(y,\cdot)dy}=0,}{}
\end{array}
\right.
\end{equation}
So the only difference in the definitions of $\bw^{(\bp)}$ and $\bv^{(\bp)}$ is that while for the first we use the matrix $A=\{a_{ij\alpha\beta}\}$,
for the latter we use the transpose (with respect to action on $\R^{md}$) matrix $A^{T}=\{a_{ji\beta\alpha}\}$.

Let $\bv^{(\bp)}_{\epsilon}(x,\omega)= \epsilon \bv^{(\bp)}( \frac{x}{\epsilon},\omega)$. Applying Lemma \ref{lem:2} to the (stationary) function
\[
g(x,\omega)= (\nabla\Phi)^{-1}(x,\omega)(\nabla\tilde\bv^{(\bp)})(x,\omega)
\]
we obtain
\begin{eqnarray}
\nabla \bv^{(\bp)}_{\epsilon}(x,\omega)&=& g(\Phi^{-1}(\frac{x}{\epsilon},\omega),\omega) \nonumber\\
&\rightharpoonup&  c_{\Phi}\E{\int_{\Phi(Q,\cdot)}g(\Phi^{-1}(x,\cdot),\cdot)dx}  \nonumber\\
&=&c_{\Phi} \E{\int_{\Phi(Q,\cdot)}\nabla\bv^{(\bp)}(x,\cdot)dx} \nonumber\\
&=& 0, \label{p3}
\end{eqnarray}
$*$-weakly in $L^{\infty}(\R^d)$ almost surely. Now, let us define
\[
\bq(y,\omega)=A^T(\Phi^{-1}(y,\omega)) (\nabla\bv^{(\bp)}(y,\omega) +\bp) \; , \quad \bq^{\epsilon}(x,\omega)=\bq(\frac{x}{\epsilon},\omega).
\]
Applying Lemma \ref{lem:2} to the function
\[
g(x,\omega)= A^T(x)\Big[    (\nabla\Phi)^{-1}(x,\omega) \nabla\tilde\bv^{(\bp)}(x,\omega) +\bp  \Big]
\]
we obtain
\begin{equation}
\bq^{\epsilon}(x,\omega) \rightharpoonup \bq^* :=c_{\Phi}\E{\int_{\Phi(Q,\cdot)}A^T(\Phi^{-1}(y)) [\nabla\bv^{(\bp)}(y,\cdot) +\bp]dy}
\label{qstar}
\end{equation}
*-weakly in $L^{\infty}(\R^d)$, almost surely. Let us compute the $j\beta$-coordinate of the first term in the last integral: for $\omega\in\Omega$
we have
\begin{eqnarray*}
\int_{\Phi(Q,\omega)}A^T(\Phi^{-1})\nabla\bv^{(\bp)}\cdot \be_{j\beta}dy &=&\int_{\Phi(Q,\omega)}A(\Phi^{-1})\be_{j\beta}
\cdot \nabla\bv^{(\bp)}dy \\
&=&- \int_{\Phi(Q,\omega)}A(\Phi^{-1})\nabla w^{j\beta}\cdot \nabla\bv^{(\bp)}dy \\
&=&\int_{\Phi(Q,\omega)}A(\Phi^{-1})\nabla w^{j\beta}\cdot \bp \, dy.
\end{eqnarray*}
Substituting in (\ref{qstar}) we obtain
\[
q^*_{j\beta} = c_{\Phi} \E{\int_{\Phi(Q,\cdot)} a_{ij\alpha\beta} +a_{jk\beta\gamma}\parder{w^{i\alpha}}{y_k}dy} p_{i\alpha},
\]
that is $\bq^*=(A^*)^T\bp$.

Since ${\diver \, \bq^{\epsilon}}$ is independent of $\epsilon$, the Compensated Compactness Theorem \cite[Lemma 1.1]{JKO} can be applied and
(\ref{p1}) together with (\ref{qstar}) imply
\begin{equation}
 \nabla \bu^{\epsilon} \cdot \bq^{\epsilon} \to \nabla\bu^*\cdot (A^*)^T \bp \; , \quad \mbox{$*$-weakly in  $L^1(D)$ a.s.}
\label{e1}
\end{equation}
Similarly, (\ref{p2}) and (\ref{p3}) imply that
\begin{equation}
 \br^{\epsilon}\cdot (\nabla\bv^{(\bp)}_{\epsilon} +\bp) \to \br^* \cdot\bp \; , \quad \mbox{$*$-weakly in  $L^1(D)$ a.s.}
\label{e2}
\end{equation}
(We refer to \cite{JKO} for the precise definition of $*$-weak convergence in $L^1$). Since the left-hand sides of (\ref{e1}) and (\ref{e2}) are equal, we obtain that $\bu^*\cdot (A^*)^T \bp =\br^* \cdot\bp$.
This being valid for all $\bp\in\R^{md}$, we obtain that $\br^* = A^* \nabla\bu^*$, as required.$\hfill \Box$

\section{An application: Homogenization of Maxwell's equations for random bi-anisotropic media}

In this section we shall see how the results of Section \ref{sec:ell} can be applied in the study of a homogenization problem for Maxwell's equations in a random medium. The corresponding deterministic problem has been studied, e.g. in
\cite{BS,Bos,StrYan}; some preliminary results for the random case are also included in the latter, as well as in
\cite{thebook}.

In what follows we shall denote by $\hat{\phi}(p)$ the Laplace transform of a function $\phi(t)$, $t>0$. Hence the variable $p$ is complex, and typically we shall have $p\in\C^+ :=\{ {\rm Re}\, p>0\}$.

\subsection{The Maxwell system}

In a bounded domain $\cO\subset\R^3$ with Lipschitz boundary we consider for a fixed $\epsilon>0$ the initial boundary value problem for Maxwell's
equations
\begin{eqnarray}
&& \partial_t \bD^{\epsilon} =\curl \bH^{\epsilon} +\bF(x,t)\nonumber \\
&& \partial_t \bB^{\epsilon} =-\curl \bE^{\epsilon}+\bG(x,t) ,\quad\qquad  x\in \cO , \; t>0, \label{eq:30}\\
&& \bE^{\epsilon}( x,0)=\bz,\quad \bH^{\epsilon}( x,0)=\bz ,
\qquad\qquad  x\in \cO, \nonumber \\
&& \bn\times \bE^{\ep}=\bz ,\qquad\qquad
\qquad\qquad\quad  x\in\partial \cO ,\; t>0.
\nonumber
\end{eqnarray}
The system (\ref{eq:30}) is accompanied by linear constitutive laws of the form
\begin{eqnarray}
&&\bD^{\epsilon}=\varepsilon^{\ep}_0\bE^{\ep}+ \xi_0^{\ep}\bH^{\ep}+\varepsilon_{d}^{\ep}* \bE^{\ep}
+\xi_{d}^{\ep}*\bH^{\ep}\nonumber\\
&& \bB^{\epsilon}=\zeta_0^{\ep}\bE^{\ep} +\mu_0^{\ep}\bH^{\ep}+\zeta_{d}^{\ep}*\bE^{\ep}
+\mu_{d}^{\ep}*\bH^{\ep},
\label{eq:35}
\end{eqnarray}
describing the anisotropic media; here and below the symbol $*$ stands for temporal convolution. In compact notation the constitutive laws can be written as
\[
\cD^{\epsilon} =A_0^{\epsilon}\cE^{\epsilon} + A^{\epsilon}_d * \cE^{\epsilon},
\]
where $\cD^{\epsilon} =[ \bD^{\epsilon} , \bB^{\epsilon}]^T$, $\cE^{\epsilon} =[ \bE^{\epsilon} , \bH^{\epsilon}]^T$ and
\[
A_0^{\epsilon}(x)=\left[
\begin{array}{cc}
\varepsilon_0^{\ep}(x)  & \xi_0^{\ep}(x) \\
\zeta_0^{\ep}(x)  & \mu_0^{\ep}(x)
\end{array}
\right]
 \quad ,\qquad
A^{\epsilon}_d(x,t)=\left[
\begin{array}{cc}
\varepsilon^{\ep}_d(x,t)  & \xi^{\ep}_d(x,t) \\
\zeta^{\ep }_d(x,t)  & \mu^{\ep}_d(x,t)
\end{array}
\right].
\]
We assume that $A_0^{\epsilon}\in (L^{\infty}(\cO))^{6\times 6}$ and $A^{\epsilon}_d\in (L^{\infty}(\cO\times (0,\infty)))^{6\times 6}$
and that the following ellipticity conditions are satisfied:
\begin{equation}
A_0^{\epsilon}(x)\cU\cdot \cU \geq c|U|^2  \;\; , \qquad  A^{\epsilon}_d(x,t)\cU\cdot \cU \geq 0  \; , \quad
\mbox{ for all $\cU\in\R^6$},
\label{ellipticity}
\end{equation}
for some $c>0$ and all $x\in\cO$, $t>0$ and $\epsilon>0$

\subsection{Homogenization of the random Maxwell system}

We now make the assumption that the matrices $A_0^{\epsilon}(x)$ and $A^{\epsilon}_d(x,t)$ above are random, of the form studied in Section
\ref{sec:ell}, that is
\[
A_0^{\epsilon}(x,\omega) = A_0(\Phi^{-1}(\frac{x}{\epsilon},\omega)) \quad , \qquad
A_d^{\epsilon}(x,t,\omega) = A_d(\Phi^{-1}(\frac{x}{\epsilon},\omega),t),
\]
where $\Phi$ is a random mapping which is a diffeomorphism almost surely, with stationary gradient and
$A_0(y)$ and $A_d(y,t)$ are periodic in $y\in\R^3$ with period cell $Q$ and satisfy
\begin{equation}
A_0(y)\cU\cdot \cU \geq c|U|^2  \;\; , \qquad  A_d(y,t)\cU\cdot \cU \geq 0  \; , \quad
\mbox{ for all $\cU\in\R^6$},
\label{ellipticity1}
\end{equation}
for some $c>0$ and all $y\in Q$ and $t>0$. We note that it follows from (\ref{ellipticity1}) that the matrix
\begin{equation}
\tilde{A}(y,p) := A_0(y)+ \hat{A}_d(y,p) =
\left[ \begin{array}{cc} \varepsilon_0(y) +\hat{\varepsilon}_d(y,p) & \xi_0(y)+\hat{\xi}_d(y,p)  \\
\zeta_0(y) +\hat{\zeta}_d(y,p) & \mu_0(y) +\hat{\mu}_d(y,p) \end{array} \right] =:
\left[ \begin{array}{cc} \tilde{\varepsilon}(y,p) & \tilde{\xi}(y,p)  \\
{\tilde{\zeta}}(y,p)  & \tilde{\mu}(y,p) \end{array} \right]
\label{eq:100}
\end{equation}
satisfies
\[
\inprod{\tilde{A}(y,p)\cU}{\cU} \geq c |\cU|^2,\quad
\by\in \R^3, \; p\in\C_+ ,\;\cU\in\R^{6}.
\]
We make the following assumption:

{\bf Assumption 1.}
The Maxwell system (\ref{eq:30}) - (\ref{eq:35}) is uniquely solvable for all $\epsilon>0$ and $\omega\in\Omega$ and the solution
satisfies $\|\bE^{\ep}\|_{L^2(\cO)},\|\bH^{\ep}\|_{L^2(\cO)}, \|\bD^{\ep}\|_{L^2(\cO)},\|\bB^{\ep}\|_{L^2(\cO)}\leq c$ for all $\epsilon,t>0$ and $\omega\in\Omega$.

For a variety of natural conditions under which Assumption 1 is valid we refer to \cite{thebook}.
E.g., if we only assume that $\|\bE^{\ep}\|_{L^2(\cO)},\|\bH^{\ep}\|_{L^2(\cO)} \leq c$ then the inequalities $\|\bD^{\ep}\|_{L^2(\cO)},\|\bB^{\ep}\|_{L^2(\cO)} \leq c$
follow if $\int_0^{\infty}\| A_d(y,t)\| dt <+\infty$, uniformly in $y\in Q$.

To state our theorem on the homogenization of the Maxwell system (\ref{eq:30}) - (\ref{eq:35}) we need to define a certain {\em homogenized} coefficient matrix $\tilde{A}^*$. This will be a 6x6 matrix depending only on $p\in\C_+$.  The matrix $\tilde{A}^*$ shall be written in block form as
\[
\tilde{A}^* =
\left[ \begin{array}{cc}
\tilde{\varepsilon}^*   &  \tilde{\xi}^* \\
\tilde{\zeta}^*  & \tilde{\mu}^*
 \end{array} \right];
\]
we define $\tilde{A}^*$ to be the transpose matrix of the limit in the sense of Theorem \ref{thm:elliptic_homog} of the sequence of transpose matrices
\[
\tilde{A}(\Phi^{-1}(\frac{x}{\epsilon},\omega),p)^{\perp} =
\left[ \begin{array}{cc} \tilde{\varepsilon}(\Phi^{-1}(\frac{x}{\epsilon},\omega),p)^{\perp} & \tilde{\zeta}(\Phi^{-1}(\frac{x}{\epsilon},\omega),p)^{\perp}  \\
{\tilde{\xi}}(\Phi^{-1}(\frac{x}{\epsilon},\omega),p)^{\perp}  & \tilde{\mu}(\Phi^{-1}(\frac{x}{\epsilon},\omega),p)^{\perp}
\end{array} \right].
\]
The precise expression for $\tilde{A}^*$ is as follows. By an application of Proposition \ref{prop:1} of Section \ref{sec:ELLHOM} there exist unique (modulo random constants) functions $u^j_1,u^j_2,v^j_1$ and $v^j_2$, $j=1,2,3$, defined by the relations
\[
\left\{
\begin{array}{l} -\parder{}{y_i}\Big\{ \tilde\varepsilon_{ik}\parder{u_1^j}{y_k} +\tilde\xi_{ik}\parder{u_2^j}{y_k}  \Big\} =\parder{\tilde\varepsilon_{ij}}{y_i} \\[0.2cm]
 -\parder{}{y_i}\bigg\{ \tilde\zeta_{ik}\parder{u_1^j}{y_k} +\tilde\mu_{ik}\parder{u_2^j}{y_k} \bigg\} =\parder{\tilde\zeta_{ij}}{y_i}
\end{array}
\right.
\]
and
\[
\left\{
\begin{array}{l} -\parder{}{y_i}\Big\{ \tilde\varepsilon_{ik}\parder{v_1^j}{y_k} +\tilde\xi_{ik}\parder{v_2^j}{y_k}  \Big\} =\parder{\tilde\xi_{ij}}{y_i} \\[0.2cm]
 -\parder{}{y_i}\bigg\{ \tilde\zeta_{ik}\parder{v_1^j}{y_k} +\tilde\mu_{ik}\parder{v_2^j}{y_k} \bigg\} =\parder{\tilde\mu_{ij}}{y_i}
\end{array}
\right.
\]
in $\cO$, where we also require that
\[
u_{\gamma}^j(y,\omega)=\tilde{u}_{\gamma}^j(\Phi^{-1}(y,\omega),\omega)  \; , \quad v_{\gamma}^j(y,\omega)=\tilde{v}_{\gamma}^j(\Phi^{-1}(y,\omega),\omega),
\]
with $\tilde{u}_{\gamma}^j$ and $\tilde{v}_{\gamma}^j$ stationary and $\E{\int_{\Phi(Q)}\nabla u^j_{\gamma}(y,\cdot)dy}=
\E{\int_{\Phi(Q)}\nabla v^j_{\gamma}(y,\cdot)dy}=0$, $\gamma=1,2$, $j=1,2,3$.
Using (\ref{eq:def_hom_m}) one can then see that
\begin{eqnarray*}
\tilde{\varepsilon}_{ij}^*&=& \det\bigg( \E{\int_Q \nabla\Phi(y,\cdot)dy}\bigg)^{-1}
\E{\int_{\Phi(Q,\cdot)} \Big(\tilde{\varepsilon}_{ij}+ \tilde{\varepsilon}_{ik} \parder{u^j_1}{y_k}
+\tilde{\xi}_{ik}\parder{u^j_2}{y_k}\Big)dy }\\[0.2cm]
\tilde{\xi}_{ij}^* &=&\det\bigg( \E{ \int_Q \nabla\Phi(y,\cdot)dy }\bigg)^{-1}
\E{\int_{\Phi(Q,\cdot)} \Big(\tilde{\xi}_{ij} +\tilde{\varepsilon}_{ik}\parder{v^j_1}{y_k} +\tilde{\xi}_{ik}\parder{v^j_2}{y_k} \Big)dy} \\[0.2cm]
\tilde{\zeta}_{ij}^* &=& \det\bigg( \E{\int_Q \nabla\Phi(y,\cdot)dy}\bigg)^{-1}
\E{ \int_{\Phi(Q,\cdot)} \Big(\tilde{\zeta}_{ij}
 + \tilde{\zeta}_{ik}\parder{u_1^j}{y_k}+\tilde{\mu}_{ik}\parder{u^j_2}{y_k} \Big)dy}\\
\tilde{\mu}_{ij}^* &=&\det\bigg( \E{\int_Q \nabla\Phi(y,\cdot)dy}\bigg)^{-1}
\E{ \int_{\Phi(Q,\cdot)} \Big(\tilde{\mu}_{ij} +\tilde{\zeta}_{ik}\parder{v^j_1}{y_k} + \tilde{\mu}_{ik}\parder{v^j_2}{y_k} \Big)dy}.
\end{eqnarray*}



{\bf Assumption 2.} (inversion of Laplace transform) There exist 3x3 matrices
$\varepsilon_0^*$, $\xi_0^*$, $\zeta_0^*$ and $\mu_0^*$ and 3x3 matrix-valued
functions $\varepsilon_d^*(t)$, $\xi_d^*(t)$, $\zeta_d^*$ and $\mu_d^*(t)$ such that
\begin{eqnarray*}
&& \varepsilon_{0,ij}^* + \hat{\varepsilon}_{d,ij}^* (p)=\tilde{\varepsilon}_{ij}^*(p) \; , \qquad
\xi_{0,ij}^* + \hat{\xi}_{d,ij}^* (p)=\tilde{\xi}_{ij}^*(p) \; , \\
&&\zeta_{0,ij}^* + \hat{\zeta}_{d,ij}^* (p)=\tilde{\zeta}_{ij}^*(p) \; , \qquad
\mu_{0,ij}^* + \hat{\mu}_{d,ij}^* (p)=\tilde{\mu}_{ij}^*(p) \; , \qquad p\in\C^+.
\end{eqnarray*}
For more information concerning this assumption we refer to \cite{widder}.

We can now state the main result of this section.
\begin{theorem}
Assume that the functions $A_0$ and $A_d(\cdot,t)$, $t>0$ belong in $C^{\delta}(\R^d)$ for some $\delta>0$.
Then the solution $[\bE^{\ep},\bH^{\ep}]$ of the Maxwell system satisfies
\[\bE^{\ep} \to \bE^*, \quad \bH^{\ep}\to \bH^*\quad\mbox{$*$-weakly in }
L^{\infty}((0,\infty)\times \Omega, L^2(\cO)),\]
where $[\bE^*,\bH^*]$ is the unique solution of the Maxwell system
\begin{eqnarray}
&& \partial_t \bD^* =\curl \bH^* +\bF(x,t) \nonumber \\
&& \partial_t \bB^* =-\curl \bE^* +\bG(x,t) ,\qquad x\in \cO , \; t>0,
\label{eq:3*} \\
&& \bE^*( x,0)=\bz,\quad \bH^*( x,0)=\bz , \nonumber \\
&& \bn\times \bE^* =\bz,\qquad\qquad\qquad\qquad  x\in\partial \cO , \; t>0, \nonumber
\end{eqnarray}
subject to the constitutive laws
\begin{eqnarray}
&&\bD^{*}=\varepsilon^{*}_0\bE^{*}+ \xi_0^{*}\bH^{*}+\varepsilon_{d}^{*}* \bE^{*}
+\xi_{d}^{*}*\bH^{*}\nonumber\\
&& \bB^{*}=\zeta_0^{*}\bE^{*} +\mu_0^{*}\bH^{*}+\zeta_{d}^{*}*\bE^{*}
+\mu_{d}^{*}*\bH^{*}.
\label{eq:3hcl}
\end{eqnarray}
\label{thm:2}
\end{theorem}
{\em Proof.} By Assumption 1 there exist $\bE^*,\bH^*,\bD^*,\bB^*
\in L^{\infty}((0,\infty)\times \Omega, L^2(\cO))$ such that, up to
taking a subsequence $\epsilon\to 0$, there holds
\begin{equation}
\left.\begin{array}{cc}
\bE^{\epsilon}\to \bE^*, &  \bH^{\epsilon}\to \bH^* \\
\bD^{\epsilon}\to \bD^*, &  \bB^{\epsilon}\to \bB^*
\end{array} \right\} \mbox{ $*$- weakly in }L^{\infty}((0,\infty)\times \Omega, L^2(\cO)).
\label{eq:conver}
\end{equation}
We define the random 3x3 matrix-valued functions $\tilde{\varepsilon}^{\epsilon}$, $\tilde{\xi}^{\epsilon}$,
$\tilde{\zeta}^{\epsilon}$ and $\tilde{\mu}^{\epsilon}$ by
\begin{eqnarray}
&& \tilde{\varepsilon}^{\epsilon}(x,p,\omega) =\tilde{\varepsilon} (\Phi^{-1}(\frac{x}{\epsilon},\omega),p),\qquad
\tilde{\xi}^{\epsilon}(x,p,\omega) =\tilde{\xi}(\Phi^{-1}(\frac{x}{\epsilon},\omega),p),\nonumber \\
&& \tilde{\zeta}^{\epsilon}(x,p,\omega) =\tilde{\zeta}(\Phi^{-1}(\frac{x}{\epsilon},\omega),p)   , \qquad
\tilde{\mu}^{\epsilon}(x,p,\omega) =\tilde{\mu}(\Phi^{-1}(\frac{x}{\epsilon},\omega),p) .
\label{mess}
\end{eqnarray}
Taking the Laplace transform of the constitutive laws (\ref{eq:35}) with respect to the time variable and exploiting the fact that the Laplace transform turns convolutions into products, we obtain
\begin{eqnarray}
&&\hat{\bD}^{\ep} =\tilde{\varepsilon}^{\epsilon}\hat{\bE}^{\ep} +\tilde{\xi}^{\epsilon}\hat{\bH}^{\ep} \nonumber\\
&&\hat{\bB}^{\ep} =\tilde{\zeta}^{\epsilon}\hat{\bE}^{\ep} +\tilde{\xi}^{\epsilon}\hat{\bH}^{\ep},
\qquad\quad  x\in \cO\; , \; \omega\in\Omega\; , \;p\in\C_+ \; , \;
\epsilon >0.
\label{ltcl}
\end{eqnarray}

We have $\| \hat{\bE}^{\ep}\|_{L^2(\cO)}\leq \int_0^{\infty} |e^{-pt}| \|\bE^{\ep}\|_{L^2(\cO)} dt$. An analogous relation is true for $\hat{\bH}^{\ep}$, $\hat{\bD}^{\ep}$ and $\hat{\bB}^{\ep}$, hence
Assumption 1 implies
\begin{equation}
\|\hat{\bE}^{\epsilon}\|_{L^{2}(\cO)} \le C^{'}, \quad  \| \hat{\bH}^{\epsilon}\|_{L^{2}(\cO)} \le C^{'}, \quad
\|\hat{\bD}^{\epsilon}\|_{L^{2}(\cO)} \le C^{'}, \quad  \| \hat{\bB}^{\epsilon}\|_{L^{2}(\cO)} \le C^{'},
\label{blt}
\end{equation}
where $C'$ is independent of $\omega\in\Omega$ and $\epsilon>0$ (but not of $p\in\C^+$).
It then easily follows that for fixed $p\in\C_+$
\begin{equation}
\left.\begin{array}{cc}
\hat{\bE}^{\epsilon}\wsc \hat{\bE}^* \; , &  \hat{\bH}^{\epsilon}\wsc \hat{\bH}^* \\
\hat{\bD}^{\epsilon}\wsc \hat{\bD}^* \; , &  \hat{\bB}^{\epsilon}\wsc \hat{\bB}^*
\end{array} \right\} \mbox{$*$ - weakly in }(L^{\infty}(\Omega,L^2(\cO)))^3 \quad
\label{eq:conmax}
\end{equation}
but also weakly in $L^2(\cO)$ almost surely.
Now, taking the Laplace transform of Maxwell's equations (\ref{eq:30}) we obtain
\begin{eqnarray}
&& p \, \hat{\bD}^{\epsilon} =\curl \hat{\bH}^{\epsilon} + \hat{\bF}, \nonumber \\
&& p \, \hat{\bB}^{\epsilon} = -\curl \hat{\bE}^{\epsilon} + \hat{\bG};
\label{ch1}
\end{eqnarray}
hence (\ref{blt}) implies that the $L^{2}(\cO)$-norms of $\curl\hat{\bE}^{\epsilon}$, $\curl\hat{\bH}^{\epsilon}$ are bounded and employing the same weak compactness argument we have the existence of $\phi, \;\psi$ such that
 \begin{eqnarray}
\curl \hat{\bE}^{\epsilon} \wsc \phi,
 \nonumber \\
\curl \hat{\bH}^{\epsilon} \wsc \psi,
 \nonumber
 \end{eqnarray}
 (up to subsequences) $*$-weakly in $(L^{\infty}(\Omega , L^{2}(\cO)))^3$. Standard arguments based on uniqueness of weak limits allow us to identify $\phi=\curl E^{*}$ and $\psi =\curl H^{*}$. The above considerations lead to the conclusion that for each $p\in\C^+$,
 \begin{eqnarray}
&& \hat{\bE}^{\epsilon} \rightharpoonup \hat{\bE}^{*}, \nonumber \\
&& \hat{\bH}^{\epsilon}\rightharpoonup \hat{\bH}^{*}
 \label{ch3}
 \end{eqnarray}
$*$-weakly in $L^{\infty}(\Omega , H(\curl , \cO))$, and also in $H(\curl , \cO)$ almost surely.
We now take the limit $\epsilon \rightarrow 0$ in (\ref{ch1})
(weakly in $(L^2(\cO))^3$, for fixed $p\in\C^+$ and $\omega\in\Omega$) and obtain that
 \begin{eqnarray}
 p \, \hat{\bD}^{*} &=&\curl \hat{\bH}^{*} + \hat{\bF}^*,
 \nonumber \\
 p \, \hat{\bB}^{*} &=& -\curl \hat{\bE}^{*} + \hat{\bG}^*.
 \label{ch2}
 \end{eqnarray}
This implies that $\bE^*,\bH^*,\bD^*$ and $\bB^*$ are solutions of the Maxwell system
\begin{eqnarray}
&& \partial_t \bD^* =\curl \bH^* +\bF^*(t) \nonumber \\
&& \partial_t \bB^* =-\curl \bE^* +\bG^*(t) ,\qquad\qquad x\in \cO ,\;  t>0,
\label{eq:2h} \\
&& \bE^*(x,0)=\bz,\quad \bH^*(x,0)=\bz, \qquad\quad x\in \cO.
\end{eqnarray}
Hence it remains to establish that the boundary
condition $\bn\times\bE^*=\bz$ is satisfied
and that the vector fields $\bE^*,\bH^*, \bD^*$ and $\bB^*$ are related by
the constitutive laws (\ref{eq:3hcl}).

{\em Validity of the boundary condition.} We first note that the
boundary condition is understood in the sense of the trace
operator $H(\curl, \cO)\to H^{-\frac{1}{2}} (\partial \cO)$,
$\bU\mapsto \bn\times \bU |_{\partial \cO}$. Let us fix a
function ${\phib}\in H^{\half}(\partial \cO)$. There exists
\cite[p. 341]{DL} ${\bf\Phi}\in H^1( \cO)$ such that
${\bf\Phi}|_{\partial \cO}=\phib$. Now, for $\ep>0$ there holds
\begin{eqnarray*}
&& \int_{\cO}\curl{\bf\Phi}\cdot
\bE^{\ep}=\int_{\cO}\curl\bE^{\ep}\cdot{\bf\Phi}
+\int_{\partial \cO}{\bf\Phi} (\bn\times\bE^{\ep}) , \\
&& \int_{\cO}\curl{\bf\Phi}\cdot
\bE^*=\int_{\cO}\curl\bE^*\cdot{\bf\Phi}
+\int_{\partial \cO}{\bf\Phi} (\bn\times\bE^{*}).
\end{eqnarray*}
Combining these with the fact that
$\bn\times\bE^{\ep}|_{\partial \cO}=\bz$ and using the relations
\begin{eqnarray*}
&& \int_{\cO}\curl{\bf\Phi}\cdot\bE^{\ep} \to \int_{\cO}\curl{\bf\Phi}\cdot\bE^* \\
&& \int_{\cO}\curl\bE^{\ep}\cdot{\bf\Phi} \to
\int_{\cO}\curl\bE^*\cdot{\bf\Phi}, \qquad\quad (\ep\to 0)
\end{eqnarray*}
we obtain
\[\int_{\partial \cO}\phib (\bn\times\bE^{*})=\int_{\partial \cO}{\bf\Phi} (\bn\times\bE^{*})
=0.\] Since $\phib\in H^{\half}(\partial \cO)$ was arbitrary, we
conclude that $\bn\times\bE^{*} =\bz$ on $\partial \cO$.

{\em Validity of the constitutive laws.} In order to prove the validity of the constitutive laws we shall need to consider an auxiliary elliptic system and apply Theorem \ref{thm:elliptic_homog}. The auxiliary system will have as a coefficient matrix the transpose of the matrix (cf. (\ref{mess}))
\[
\left[
\begin{array}{cc}
 \tilde{\varepsilon}^{\epsilon} & \tilde{\xi}^{\epsilon} \\
 \tilde{\zeta}^{\epsilon}  & \tilde{\mu}^{\epsilon}
\end{array}
\right] ;
\]
the system will also depend on the parameter $p\in\C^+$, as well as on $\omega\in\Omega$.

Specifically, let us fix a domain $V$ with smooth boundary compactly contained in $\cO$. We recall definition (\ref{eq:100})
of the $3\times 3$ periodic matrix-valued functions $\tilde{\varepsilon}$, $\tilde{\xi}$, $\tilde\zeta$ and $\tilde{\mu}$  and
we define the random elliptic operator $L^{\epsilon} :  (H_{0}^{1}(V))^2 \rightarrow (H^{-1}(V))^2$ by
 \begin{eqnarray}
 L^{\epsilon}=\left[
 \begin{array}{cccc}
 -\diver( (\tilde{\varepsilon}^{\epsilon})^{\perp} \, \grad) & -\diver((\tilde{\zeta}^{\epsilon})^{\perp} \, \grad)
 \nonumber \\
 -\diver((\tilde{\xi}^{\epsilon})^{\perp} \, \grad) & -\diver((\tilde{\mu}^{\epsilon})^{\perp} \, \grad)
 \nonumber
 \end{array}
 \right]
 \end{eqnarray}
where the $3\times 3$ random matrix-values functions $\tilde{\varepsilon}^{\epsilon}$, $\tilde{\xi}^{\epsilon}$,
$\tilde{\zeta}^{\epsilon}$ and $\tilde{\mu}^{\epsilon}$ have been defined by (\ref{mess}). Here and below we denote by
${\grad }$ the usual gradient operator in $\R^3$.
Let $L^*: (H^1_0(V))^2\to (H^{-1}(V))^2$ be the operator given by
\[
L^* =\left[ \begin{array}{cc}
-\diver((\tilde{\varepsilon}^*)^{\perp}\grad)  & -\diver((\tilde{\zeta}^*)^{\perp}\grad) \\
-\diver((\tilde{\xi}^{*})^{\perp}\grad)  &  -\diver((\tilde{\mu}^*)^{\perp}\grad)
\end{array}\right].
\]
We note that the coefficients of $L^*$ depend only on the parameter $p\in\C_+$.

Let $\cG =[g_1,g_2]^T
\in (L^2(V))^2$ be fixed. By Theorem \ref{thm:elliptic_homog} and the definition of $\tilde{A}^*$ the solutions
${\cal U}^{\epsilon}=[u^{\epsilon} , v^{\epsilon}]^{T}$, ${\cal U}^{*}=[u^{*} , v^{*}]^{T}$ of
\begin{equation}
L^{\epsilon} {\cal U}^{\epsilon}= {\cal G}, \qquad\qquad
L^{*} {\cal U}^{*}= {\cal G},
\label{ggg}
\end{equation}
satisfy
 \begin{eqnarray}
 \grad \, u^{\epsilon}  \rightharpoonup \grad \, u^{*},
 \nonumber \\
 \grad \, v^{\epsilon}  \rightharpoonup \grad \, v^{*},
 \label{ym}
 \end{eqnarray}
weakly in $L^{2}(V)$ a.s., and also
 \begin{eqnarray}
 (\tilde{\varepsilon}^{\epsilon})^{\perp} \,  \grad \, u^{\epsilon} + (\tilde{\zeta}^{\epsilon})^{\perp}  \grad \, v^{\epsilon} \rightharpoonup (\tilde{\varepsilon}^{*})^{\perp} \, \grad \, u^{*} + (\tilde{\zeta}^{*})^{\perp} \grad \, v^{*} ,
 \nonumber \\
( \tilde{\xi}^{\epsilon})^{\perp} \grad \, u^{\epsilon}  + (\tilde{\mu}^{\epsilon})^{\perp}  \grad \, v^{\epsilon}\rightharpoonup
(\tilde{\xi}^{*})^{\perp} \, \grad \, u^{*} + (\tilde{\mu}^{*})^{\perp} \grad \, v^{*} ,
 \label{eqn:CONV2}
 \end{eqnarray}
 weakly in $L^{2}(V)$ a.s. .
 The vector identity $\curl \grad =0$ together with (\ref{ym}) imply that in fact
 \begin{eqnarray}
 \grad \, u^{\epsilon}  \rightharpoonup \grad \, u^{*},
 \nonumber \\
 \grad \, v^{\epsilon}  \rightharpoonup \grad \, v^{*},
 \label{eqn:CONVB}
 \end{eqnarray}
 weakly in $H(\curl , V)$, a.s.. Moreover, the identity $\diver\curl =0$ together with (\ref{ch1}) and (\ref{ch2}) yields
$\diver \hat{\bD}^{\epsilon}=\diver \hat{\bD}^*$ and $\diver \hat{\bB}^{\epsilon}=\diver \hat{\bB}^*$; this combined with (\ref{eq:conmax}) implies
 \begin{eqnarray}
 \hat{\bD}^{\epsilon}  \rightharpoonup \hat{\bD}^{*},
 \nonumber \\
\hat{\bB}^{\epsilon}  \rightharpoonup \hat{\bB}^{*},
 \label{eqn:div}
 \end{eqnarray}
 weakly in $H(V,\diver)$, a.s.. Relations (\ref{eqn:CONVB}), (\ref{eqn:div}) and the div-curl lemma now imply that
\begin{eqnarray}
 \hat{\bD}^{\epsilon} \cdot \grad u^{\epsilon} \rightarrow \hat{\bD}^{*} \cdot \grad u^{*},
 \nonumber \\
 \hat{\bB}^{\epsilon} \cdot \grad v^{\epsilon} \rightarrow \hat{\bB}^{*} \cdot \grad v^{*},
 \label{ch4}
 \end{eqnarray}
 a.s. in ${\cal D}^{'}(V)$. Adding up (\ref{ch4}) we obtain
 \begin{eqnarray}
 \hat{\bD}^{\epsilon} \cdot \grad \, u^{\epsilon}  + \hat{\bB}^{\epsilon} \cdot \ \grad \, v^{\epsilon}
 \rightarrow \hat{\bD}^{*} \cdot \grad \, u^{*}  + \hat{\bB}^{*} \cdot \cdot \grad \, v^{*}
 \label{eqn:CONVE}
 \end{eqnarray}
a.s. in ${\cal D}^{'}(V)$.
Moreover we have
\[
\begin{array}{c}
 -\diver \Big((\tilde{\varepsilon}^{\epsilon})^{\perp}\grad u^{\epsilon} +(\tilde{\zeta}^{\epsilon})^{\perp}\grad v^{\epsilon} \Big)= g_1
=-\diver \Big((\tilde{\varepsilon}^*)^{\perp} \grad u^* +(\tilde{\zeta}^*)^{\perp}\grad v^*\Big) \\[0.2cm]
-\diver \Big((\tilde{\xi}^{\epsilon })^{\perp}\grad u^{\epsilon} +(\tilde{\mu}^{\epsilon})^{\perp}\grad v^{\epsilon} \Big)= g_2
=-\diver \Big((\tilde{\zeta}^{*})^{\perp} \grad u^* +(\tilde{\mu}^*)^{\perp}\grad v^* \Big);
\end{array}
\]
these together with (\ref{eqn:CONV2}) imply
\[
\left. \begin{array}{c}
(\tilde{\varepsilon}^{\epsilon})^{\perp}\grad u^{\epsilon} +(\tilde{\zeta}^{\epsilon})^{\perp}\grad v^{\epsilon}
\rightharpoonup (\tilde{\varepsilon}^*)^{\perp}\grad u^*+(\tilde{\zeta}^*)^{\perp}\grad v^* \\[0.2cm]
(\tilde{\xi}^{\epsilon})^{\perp}\grad u^{\epsilon} +(\tilde{\mu}^{\epsilon})^{\perp}\grad v^{\epsilon}
\rightharpoonup (\tilde{\zeta}^{*})^{\perp}\grad u^*+(\tilde{\mu}^*)^{\perp}\grad v^*
\end{array}\right.
\]
weakly in $H(V,\diver)$, almost surely. Combining this with (\ref{ch3}) we obtain by another application of the div-curl lemma that
 \begin{eqnarray}
&& \big((\tilde{\varepsilon}^{\epsilon})^{\perp} \,  \grad \, u^{\epsilon} + (\tilde{\zeta}^{\epsilon})^{\perp}  \grad \, v^{\epsilon}\big) \cdot \hat{\bE}^{\epsilon} \rightarrow \big( (\tilde{\varepsilon}^{*})^{\perp} \grad \, u^{*} + (\tilde{\zeta}^{*})^{\perp} \grad \, v^{*} \big) \cdot \hat{\bE}^{*} ,
 \nonumber \\[0.2cm]
&& \big( (\tilde{\xi}^{\epsilon})^{\perp} \grad \, u^{\epsilon}  + (\tilde{\mu}^{\epsilon})^{\perp}  \grad \, v^{\epsilon}\big) \cdot
 \hat{\bH}^{\epsilon} \rightarrow
 \big( (\tilde{\xi}^{*})^{\perp} \, \grad \, u^{*} + (\tilde{\mu}^{*})^{\perp} \grad \, v^{*}  \big) \cdot \hat{\bH}^{*} ,
 \label{eqn:CONV3}
 \end{eqnarray}
almost surely in ${\cal D}^{'}(V)$.
Adding relations (\ref{eqn:CONV3}) we obtain
\begin{eqnarray}
&&  (\tilde{\varepsilon}^{\epsilon} \hat{\bE}^{\epsilon} +\tilde{\xi}^{\epsilon }\hat{\bH}^{\epsilon} ) \cdot \grad \, u^{\epsilon} +
  (\tilde{\zeta}^{\epsilon} \hat{\bE}^{\epsilon}+\tilde{\mu}^{\epsilon} \hat{\bH}^{\epsilon} ) \cdot \grad \, v^{\epsilon}
   \rightarrow   \label{eqn:HOMCONV} \\
 && \hspace{3.5cm}(\tilde{\varepsilon}^{*} \hat{\bE}^{*} +\tilde{\xi}^{*}
  \hat{\bH}^{*} ) \cdot \grad \, u^{*} +
  (\tilde{\zeta}^{*} \hat{\bE}^{*}+\tilde{\mu}^* \hat{\bH}^{*} ) \cdot \grad \, v^{*}
  \nonumber
  \end{eqnarray}
in ${\cal D}^{'}(V)$, almost surely, or equivalently by (\ref{ltcl}),
  \begin{eqnarray}
  \hat{\bD}^{\epsilon} \cdot \grad \, u^{\epsilon}  + \hat{\bB}^{\epsilon} \cdot \grad \, v^{\epsilon}
   \rightarrow
  (\tilde{\varepsilon}^{*} \hat{\bE}^{*} +\tilde{\xi}^{*}
  \hat{\bH}^{*} ) \cdot \grad \, u^{*} +
  (\tilde{\xi}^{* T} \hat{\bE}^{*}+\tilde{\mu}^{*} \hat{\bH}^{*} ) \cdot \grad \, v^{*}
  \label{eqn:CONVC}
  \end{eqnarray}
  a.s. in ${\cal D}^{'}(V)$.  Combining (\ref{eqn:CONVC}) and (\ref{eqn:CONVE}) we obtain by the uniqueness of limits,
\[
  \hat{\bD}^{*} \cdot \grad u^{*} + \hat{\bB}^{*} \cdot \grad v^* = (\tilde{\varepsilon}^{*} \hat{\bE}^{*} +\tilde{\xi}^{*}
  \hat{\bH}^{*} ) \cdot \grad \, u^{*} +
  (\tilde{\xi}^{* T} \hat{\bE}^{*}+\tilde{\mu}^{*} \hat{\bH}^{*} ) \cdot \grad \, v^{*} \; , \quad \mbox{ in }V.
\]
We claim that
\begin{eqnarray}
&&  \hat{\bD}^{*}=\tilde{\varepsilon}^{*} \,\hat{\bE}^{*} +\tilde{\xi}^{*} \hat{\bH}^{*}, \nonumber \\
 && \hat{\bB}^{*}= \tilde{\xi}^{* T} \,\hat{\bE}^{*} + \tilde{\mu}^{*} \, \hat{\bH}^{*}, \qquad \mbox{ in $V$, for
 all $p\in\C^+$.}
\label{lapl}
\end{eqnarray}
To see this, let $B(x_0,\rho)$ be a ball compactly contained in $V$ and let ${\bf e}_i\in\R^3$ be a fixed basis vector. Let $u^*\in\ C^{\infty}_0(V)$ be such that $u^*=x_i$ in $B(x_0,\rho)$ and let $v^*$ be identically zero in $V$.
Denoting as above ${\cal U}^{*}=[u^{*} , v^{*}]^{T}$, we define the vector field
${\cal U}^{\epsilon}=[u^{\epsilon} , v^{\epsilon}]^{T}\in (H^1_0(V))^2$ by requiring that
\[
L^{\epsilon} {\cal U}^{\epsilon}=  L^{*} {\cal U}^{*}  , \qquad\quad \mbox{ in }V.
\]
Hence we obtain that $\hat{\bD}^{*} \cdot \be_i  = (\tilde{\varepsilon}^{*} \hat{\bE}^{*} +\tilde{\xi}^{*}\hat{\bH}^{*} ) \cdot \be_i$ in
$B(x_0,\rho)$. Since $\be_i$ is an arbitrary basis vector and $B(x_0,\rho)$ is also arbitrary, we conclude that
$\hat{\bD}^{*}  = \tilde{\varepsilon}^{*} \hat{\bE}^{*} +\tilde{\xi}^{*}\hat{\bH}^{*}$ in $V$. Similarly we obtain the second relation in
(\ref{lapl}).

We finally note that (\ref{lapl}) is the Laplace transform of the stated constitutive laws (\ref{eq:3hcl});
since $V$ was arbitrary, this concludes the proof. \hfill $\Box$

{\bf Acknowledgment.} We thank Professor Claude Le Bris for useful discussions.


\end{document}